\documentstyle[12pt]{article}

\begin{document}

\begin{center}
{\large {\bf {Twist deformations for Yangians}}} \vskip
5mm {\sc Vladimir D. Lyakhovsky
\footnote{E-mail address: lyakhovs@snoopy.phys.spbu.ru } }\\ \vskip 2mm
Theoretical Department, St. Petersburg State University, 198904, St.
Petersburg, Russia \\ \vskip 2mm
\end{center}

\begin{abstract}
It is demonstrated how chains of twists for classical Lie algebras induce
the new twist deformations (the deformed Yangians) that quantize the
generalized rational solutions of the classical Yang-Baxter equation.
For the case of $\left\{ Y(g)| g=so(2N+1) \right\}$ the explicit expression
for the corresponding $R_{\cal Y}$-matrix in the defining representation
is given.
\end{abstract}


\section{Introduction}

Yangians $Y(g)$ were introduced by V.G.Drinfeld \cite{D3} as the
quantizations of rational solutions of the classical Yang-Baxter equation
(CYBE). It was demonstrated by A.Stolin \cite{STO1} that rational solutions
with values in the tensor square of the Lie algebra $g$ space can be written
in the form
$$
\frac{C^2}{u-v}+a_0+b_1u+b_2v+cuv
$$
with $a_0,b_i,c\in g$ and $C^2$ -- the second Casimir of the algebra $g$.
Notice that this does not necessarily signify the simplicity of $g$. The
existence of Yangians for nonsemisimple Lie algebras was predicted in \cite
{STO2} and explicitly demonstrated in \cite{L-A}.

The important class of rational solutions is of the form
$$
\frac{C^2}{u-v}+r_0,
$$
where the additional summand is a constant solution of the CYBE. As it was
proved in \cite{KULSTO} the quantization of such solutions leads to the
twist deformations $Y_{{\cal F}}(g)$ of the Yangians $Y(g)$ (the latter
refer to the canonical rational solutions $\frac{C^2}{u-v}$). The deformed
Yangian $Y_{{\cal F}}(g)$ has the same multiplication as in $Y(g)$, the
twisted coproduct
$$
\Delta _{{\cal F}}\left( y\right) ={\cal F}\Delta \left( y\right) {\cal
F}^{-1},
$$
and the transformed ${\cal R}$-matrix
$$
{\cal R_{F}}=\left( {\cal F}\right) _{21}{\cal RF}^{-1}.
$$

The twisting element ${\cal F}$ has to satisfy the equations \cite{D2}:
\begin{equation}
\begin{array}{l}
\label{TE}{\cal F}_{12}(\Delta \otimes id)({\cal F})={\cal F}_{23}(id\otimes
\Delta )({\cal F}),\\
(\epsilon \otimes id)({\cal F})=(id\otimes \epsilon )({\cal F})=1.
\end{array}
\end{equation}

When ${\cal A}$ and ${\cal B}$ are the universal enveloping algebras: ${\cal
A}=U({\mbox{g }})\supset {\cal B}=U({\mbox{l }})$ with ${\mbox{g }}\supset {
\mbox{l }}$and $U({\mbox{l }})$ is the minimal subalgebra on which ${\cal F}$
is completely defined as ${\cal F}\in U({\mbox{l }})\otimes U({\mbox{l}})$
then ${\mbox{l }}$ is called the {\bf carrier} algebra for ${\cal F}$.

The existence of a twist can be formulated in terms of a nondegenerate
bilinear form on the carrier algebra. The generators dual to the PBW basic
elements of $U(g)$ are very important. They provide the explicit
presentation of the twisting element \cite{KLM}.

The two known examples of twists that were written explicitly
correspond to the two-dimensional carrier subalgebra $B(2)$ with the
generators $H$ and $E$,
$$
[H,E]=E,
$$
and the four-dimensional carrier subalgebra ${\bf L}$:
\begin{equation}
\label{l-def}
\begin{array}{l}
[H,E]=E,\quad [H,A]=\alpha A,\quad [H,B]=\beta B, \\
[0.2cm][A,B]=E,\quad [E,A]=[E,B]=0,\quad \quad \alpha +\beta =1.
\end{array}
\end{equation}
The first one is called the {\bf Jordanian twist} \cite{OGIEV} and has the
twisting element
\begin{equation}
\label{ogiev}\Phi _{{\cal J}}=e^{H\otimes \sigma },\quad \quad \sigma =\ln
(1+E).
\end{equation}
The second one is the {\bf extended Jordanian twist} suggested in \cite{KLM}
. The corresponding twisting element
\begin{equation}
\label{t-ext}{\cal F}_{{\cal E}(\alpha ,\beta )}=\Phi _{{\cal E}(\alpha
,\beta )}\Phi _{{\cal J}}
\end{equation}
contains the {\bf Jordanian factor }$\Phi_{{\cal J}}$ and the {\bf extension}
\begin{equation}
\label{fractions}\Phi _{{\cal E}(\alpha ,\beta )}=\exp \{A\otimes Be^{-\beta
\sigma }\}.
\end{equation}

In general the composition of two twists is not a twist. But there are some
important examples of the opposite behaviour. When ${\bf L}$ is a subalgebra
${\bf L\subset }{\mbox{g }}$ there may exist several pairs of generators of
the type $(A,B)$ arranged so that the Jordanian twist can acquire several
similar extensions \cite{KLM}. This demonstrates that some twistings can be
applied successively to the initial Hopf algebra even in the case when their
carrier subalgebras are nontrivially linked. In the universal enveloping
algebras for classical Lie algebras there exists the possibility to
construct systematically the special sequences of twists called {\bf chains}
\cite{KLO}:
\begin{equation}
\label{mainchain}{\cal F}_{{\cal B}_{p\prec 0}}\equiv {\cal F}_{{\cal B}_p}
{\cal F}_{{\cal B}_{p-1}}\ldots {\cal F}_{{\cal B}_0}.
\end{equation}
The factors ${\cal F}_{{\cal B}_k}=\Phi _{{\cal E}_k}\Phi _{{\cal J}_k}$ of
the chain are the twisting elements of the extended Jordanian twists for the
initial Hopf algebra ${\cal A}_0$. Here the extensions $\left\{ \Phi _{{\cal
E}_k},k=0,\ldots p-1\right\} $ contain the fixed set of normalized factors $
\Phi _{{\cal E}(\alpha ,\beta )}=\exp \{A\otimes Be^{-\beta \sigma }\}$ ,
the {\bf full set}. It was proved that in the classical Lie
algebras that conserve symmetric invariant forms such chains can be made
maximal and proper. This means that for the algebras $U(A_N)$, $U(B_N)$ and $
U(D_N)$ of the three classical series there exist chains ${\cal F}_{{\cal B}
_{p\prec 0}}$ that cannot be reduced to a chain for a simple subalgebra and
their full sets of extensions are the maximal sets in the sense described
below.

To construct a maximal proper chain for ${\cal A}=U(g)$ (where $g$ is a
classical Lie algebra with the root system $\Lambda _{{\cal A}}$) the
sequences ${\cal A}\equiv {\cal A}_0\supset {\cal A}_1\supset \ldots \supset
{\cal A}_{p-1}\supset {\cal A}_p$ of Hopf subalgebras are to be fixed in $
{\cal A}$. For each element ${\cal A}_k$ of the sequence there must exist
the so called {\bf initial root} $\lambda _0^k$ and the set $\pi _k$ of its
{\bf constituent roots},
\begin{equation}
\label{kpi}\pi _k=\left\{ \lambda ^{\prime },\lambda ^{\prime \prime
}|\lambda ^{\prime }+\lambda ^{\prime \prime }=\lambda _0^k;\quad \lambda
^{\prime }+\lambda _0^k,\lambda ^{\prime \prime }+\lambda _0^k\neg \in
\Lambda _{{\cal A}}\right\}
\end{equation}
For any $\lambda ^{\prime }\in \pi _k$ there must be an element $
\lambda ^{\prime \prime }\in \pi _k$ that $\lambda ^{\prime }+\lambda
^{\prime \prime }=\lambda _0^k$. So, $\pi _k$ is naturally decomposed as
\begin{equation}
\pi _k=\pi _k^{\prime }\,\cup \,\pi _k^{\prime \prime },\quad \quad \pi
_k^{\prime }=\{\lambda ^{\prime }\},\quad \pi _k^{\prime \prime }=\{\lambda
^{\prime \prime }\}.
\end{equation}
The triples $\left( {\cal A}_k,\lambda _0^k,\pi _k\right) $ are subject to
the following conditions:

\begin{enumerate}
\item  $\lambda _0^k$ must be orthogonal to the roots of the subalgebra $
{\cal A}_{k-1}$,

\item  the subsets $\pi _k^{\prime }$ and $\,\pi _k^{\prime \prime }$ must
form the diagrams of conjugate representations for ${\cal A}_{k-1}$.
\end{enumerate}

In these terms the factors ${\cal F}_{{\cal B}_k}$ of the chain (\ref
{mainchain}) are fixed as follows:
\begin{equation}
\label{fabk}{\cal F}_{{\cal B}_k}=\Phi _{{\cal E}_k}\Phi _{{\cal J}_k}
\end{equation}
with
\begin{equation}
\label{jordfact}\Phi_{{\cal J}_k}=\exp \{H_{\lambda _0^k}\otimes \sigma
_0^k\},\quad \quad \sigma_0^k=\ln (1+L_{\lambda _0^k});
\end{equation}
\begin{equation}
\label{extfact}\Phi_{{\cal E}_k}=\prod_{\lambda ^{\prime }\in \pi
_k^{\prime }}\Phi_{{\cal E}_{\lambda ^{\prime }}}=\prod_{\lambda ^{\prime
}\in \pi _k^{\prime }}\exp \{L_{\lambda ^{\prime }}\otimes L_{
{\lambda _0^k}-\lambda ^{\prime }}e^{-\frac 12\sigma_0^k}\}
\end{equation}
(here $L_\lambda $ is the generator associated with the root $\lambda $).

Quantizations ${\cal A}_{{\cal F}_{{\cal B}_{p\prec 0}}}$ of classical
universal enveloping algebras produce the chains of ${\cal R_F}$-matrices:
\begin{equation}
\label{chrmat}{\cal R}_{{\cal B}_{p\prec 0}}=({\cal F}_{{\cal B}_p})_{21}(
{\cal F}_{{\cal B}_{p-1}})_{21}\ldots ({\cal F}_{{\cal B}_0})_{21}{\cal F}_{
{\cal B}_0}^{-1}\ldots {\cal F}_{{\cal B}_{p-1}}^{-1}{\cal F}_{{\cal B}
_p}^{-1}.
\end{equation}

The deformation parameters can be introduced in chains by rescaling the
generators in the subalgebras ${\cal B}_k$. Each ${\cal B}_k$ can be rescaled
separately with an independent variable $\xi _k$. When all these scaling
factors are proportional to the deformation parameter $\xi $, i.e. $\xi
_k=\xi \eta _k$, then in the classical limit the parameters $\eta _k$ appear
as the multipliers in the classical $r$-matrix:
\begin{equation}
\label{rmatr}r_{_{{\cal B}p\prec 0}}=\sum_{k=0,1,\ldots ,p}\eta _k\left(
H_{\lambda _0^k}\wedge L_{\lambda _0^k}+\sum_{\lambda ^{\prime }\in \pi
_k}L_{\lambda ^{\prime }}\wedge L_{\lambda _0^k-\lambda ^{\prime }}\right) .
\end{equation}

To obtain the necessary background for the integrable models with deformed
Yangian symmetry the simplest way is to construct the defining
representation of the universal ${\cal R_Y}$-matrix, $R=d\left( {\cal R_Y}
\right) $. In some special situations the quantum $R_{\cal F}$-matrix
corresponding to the twisted algebra can be
obtained directly from its classical counterpart $r_{\cal F}$ \cite{KULSTO}
(in particular this happens when $R=\exp \left( r\right) $ \cite{GER}). Such
simplifications are unavailable when the $r_{\cal F}$-matrices of the type
(\ref{rmatr}) are considered for the orthogonal classical Lie algebras $g$.
In such cases more information about the ${\cal R_F}$-matrix is necessary.

In this report it will be shown how the chains of extended twists (\ref
{mainchain}) lead to the series of deformed Yangians. To illustrate the
situation the explicit formulas will be presented for the case of $
g=so(2N+1) $.


\section{Chains and twisted ${\cal R}$-matrices}

Let $g$ be a classical Lie algebra of the type $B_N$ or $D_N$ and ${\cal F}_{
{\cal B}_{p\prec 0}}$ $\equiv$ ${\cal F}_{{\cal B}_p}{\cal F}_{{\cal B}
_{p-1}}\ldots $ $\ldots {\cal F}_{{\cal B}_0}$ -- its full proper chain of
extended twists.

To construct a maximal proper chain for ${\cal A}=U(g)$ the following
sequences ${\cal A}\equiv {\cal A}_0\supset {\cal A}_1\supset \ldots \supset
{\cal A}_{p-1}\supset {\cal A}_p$ of Hopf subalgebras are to be fixed:
\begin{equation}
\label{soetower}
\begin{array}{l}
U(so(2N))\supset U(so(2(N-2))\supset \ldots \supset U(so(2(N-2k))\supset
\ldots \\
\end{array}
\,\,{\rm for}\,\,D_N
\end{equation}
\begin{equation}
\label{sootower}U(so(2N+1))\supset U(so(2(N-2)+1)\supset \ldots \supset
U(so(2(N-2k)+1)\supset \ldots \,\,{\rm for}\,\,B_N
\end{equation}
In both cases the initial roots for ${\cal A}_k$ can be chosen to be $
\lambda _0^k=e_1+e_2$ (here the root subsystems are considered for ${\cal A}
_k$ separately and all the roots are written in the standard $e$-basis).

The twisting element for a chain can be written explicitly as
\begin{equation}
\label{chain}
\begin{array}{lll}
{\cal F}_{{\cal B}_{p\prec 0}} & = & \prod_{\lambda ^{\prime }\in \pi
_p^{\prime }}\left( \exp \{L_{\lambda ^{\prime }}\otimes L_{\lambda
_0^p-\lambda ^{\prime}}e^{-\frac 12\sigma _0^p}\}\right) \cdot \exp
\{H_{\lambda _0^p}\otimes \sigma _0^p\}\,\cdot \\
&  & \prod_{\lambda ^{\prime }\in \pi _{p-1}^{\prime }}\left( \exp
\{L_{\lambda ^{\prime }}\otimes L_{\lambda _0^{p-1}-
\lambda ^{\prime }}e^{-\frac12\sigma _0^{p-1}}\}\right) \cdot \exp
\{H_{\lambda _0^{p-1}}\otimes \sigma_0^{p-1}\}\,\cdot \\
&  & \ldots \\
&  & \prod_{\lambda ^{\prime }\in \pi _0^{\prime }}\left( \exp \{L_{\lambda
^{\prime }}\otimes L_{\lambda _0^0-\lambda ^{\prime }}e^{-\frac 12\sigma
_0^0}\}\right) \cdot \exp \{H_{\lambda _0^0}\otimes \sigma _0^0\}
\end{array}
\end{equation}

Let us introduce the deformation parameters $\xi _k=\xi \eta _k$ and rescale
the generators $\left\{ L_{\lambda _0^k},L_{\lambda _0^k-\lambda ^{\prime
}}\right\} $ by  $\xi _k$ in each subalgebra ${\cal A}_k$.
Using the expressions (\ref{jordfact}) and (\ref{extfact}) one can get the
first terms of the expansion for the twisting element
$$
{\cal F}_{{\cal B}_{p\prec 0}}\left( \xi \right) =I\otimes I+\xi \rho _{
{\cal B}}+{\cal O}\left( \xi ^2\right) .
$$
Here
\begin{equation}
\label{rhomat}\rho _{\cal B}=\sum_{k=0,1,\ldots ,p}\eta _k\left(
H_{\lambda _0^k}\otimes L_{\lambda _0^k}+\sum_{\lambda ^{\prime }\in \pi
_k}L_{\lambda ^{\prime }}\otimes L_{\lambda _0^k-\lambda ^{\prime }}\right)
\end{equation}
is ''a half'' of the $r_{{\cal B}_{p\prec 0}}$-matrix:
$$
r_{{\cal B}_{p\prec 0}}=\rho _{\cal B}-\tau \circ \rho _{\cal B}.
$$

The carrier subalgebras for chains of extended twists are solvable. As it
was already mentioned in Section 1 the generators of the carrier form dual
sets. In the case of extended twists one of these sets ${\cal L}=\left\{
L_{\lambda _0^k},L_{\lambda _0^k-\lambda ^{\prime }}|\lambda ^{\prime }
\in \pi_k^{\prime } \right\} $ forms a nilpotent subalgebra. In the
defining representation $d(g)$ the corresponding matrix ring is nilpotent
with the index $\kappa =3$. As it was demonstrated in \cite{KLM} the
twisting element can always be presented in the
form ${\cal F}=\exp \left( {\cal L}_i^{*}\otimes \psi^i \left( {\cal L}
\right) \right) $and the expansions for the elements $\psi^i \left(  {\cal
L}\right) $ starts with the term $\xi {\cal L}^i$. Applying these results to
the case of chains we get the following property.

{\bf Lemma.} Let $g$ be a Lie algebra of the type $A_N$, $B_N$ or $D_N$ and
$d(g)$ -- the defining representation of $g$. Then for the full proper chain
of extended twists ${\cal F}_{{\cal B}_{p\prec 0}}\left(\xi \right)=I\otimes
I+\xi \rho _{\cal B}+{\cal O}\left( \xi ^2\right)$ the following relations
are true:
\begin{enumerate}
\item  $d\left( \rho _{\cal B}^3 \right) =0,$

\item  $d\left( {\cal F}_{{\cal B}_k}\left( \xi \right) \right) =\exp \left(
\xi d\left( \rho _{{\cal B}_k} \right) \right) ,$

\item  $R_{{\cal B}_{p\prec 0}}\left( \xi \right) =d\left( {\cal R}_{{\cal B}
_{p\prec 0}}\left( \xi \right) \right) = I\otimes I-\xi d\left(
\rho _{{\cal B}}\right) +\frac{\xi ^2}2\left( d^2\left( \rho _{\cal B}
\right) +d^2\left( \tau \circ \rho _{\cal B}\right) \right) -\xi ^2d
\left( \tau \circ \rho _{\cal B}\right) d\left(\rho _{\cal B}\right).
\spadesuit $
\end{enumerate}

The construction of the $R$-matrix differs for linear and orthogonal
algebras. The most interesting is the case $g=so(M)$, i.e. the series
$B_N$ or $D_N$. Let us fix $g$ to be of $B_N$ type, $g=so(2N+1)$.
Thus ${\sf C^{2N+1}}$ is the space of the defining representation. Let
$P$ be the permutation matrix acting in ${\sf C^{2N+1}\otimes C^{2N+1}}$
$$
P=M^{\prime }\otimes M^{\prime \prime }\in {\rm Mat}\left( 2N+1,{\sf C}
\right) ^{\otimes 2}
$$
Define also the matrix
$$
K=\left( M^{\prime }\right) ^{{\sf T}}\otimes M^{\prime \prime }.
$$
It was shown in \cite{KULSTO} that if ${\cal F}$ is the twisting element for
$U(so(2N+1))$ with the twisted $R$-matrix $R_{\cal F}$ the corresponding
rational solution of the quantum Yang-Baxter equation in the defining
representation has the form
\begin{equation}
\label{yrmatr}R=d\left( {\cal R_Y}\right) =uR_{\cal F}+P-\frac
u{u+N-1/2}d\left( {\cal F}_{21}\right) Kd\left( {\cal F}^{-1}\right) .
\end{equation}
Applying the Lemma proved above we get the following final expression
for $R$.
\begin{equation}
\label{rymatr}
\begin{array}{l}
R=u\left( I\otimes I-\xi d\left( \rho _{
\cal B}\right) +\frac{\xi ^2}2d\left( \rho _{\cal B}^2+\left( \tau \circ
\rho _{\cal B}\right) ^2-2\left( \tau \circ \rho _{\cal B}\right) \rho _{
\cal B}\right) \right) +P- \\
-\frac u{u+N-1/2}d\left( I\otimes I+\tau \circ \rho _{\cal B}+\frac 12
\left( \tau \circ \rho _{\cal B}\right)^2\right) Kd\left( -I\otimes I-
\rho _{\cal B}-\frac 12\left( \rho _{\cal B}\right) ^2\right).
\end{array}
\end{equation}
This $R$-matrix describes the Yangians $Y_{{\cal F}_{{\cal B}_{p\prec 0}}}
(so(2N+1))$ deformed by the full chains of extended twists. To get the
final answer the only term that must be calculated  is the defining
representation for the $\rho _{\cal B}$-matrix (\ref{rhomat}).

\section{Example. QYBE solutions for so(2N+1)}

Let $g=so(2N+1)$. To make the illustration maximally visual we shall
consider the simplest nontrivial chain, i. e. put $p=0$ and
consider a chain with a single factor ${\cal F}_{{\cal B}_0}$ in
(\ref{mainchain}). The algorithm for other factors ${\cal F}_{{\cal B}_k}$
is similar to that of the first one and the case of a full chain can be
reconstructed using the formulas (\ref{fabk}-\ref{extfact}) and the
expressions presented below.

The matrices of the defining representation $d(g)$ are written in terms of
basic antisymmetric Okubo matrices $M_{ij}$:
$$
d\left( H_{e_i+e_j}\right) =-i\left( M_{2i-1,2i}+M_{2j-1,2j}\right) \equiv
H_{i+j},
$$
$$
d\left( L_{e_k}\right) =M_{2k,2N+1}-iM_{2k-1,2N+1}\equiv E_k,\qquad
k=1,\ldots ,N
$$
$$
\begin{array}{lcl}
d\left( L_{e_i\pm e_j}\right) &  = & \frac 12\left( -M_{2i,2j}\pm
iM_{2i,2j-1}+iM_{2i-1,2j}\pm M_{2i-1,2j-1}\right)  \\
& \equiv & E_{i+j},\qquad i<j
\end{array}
$$
In this representation the factors of the sequence
$$
d\left( {\cal F}_{{\cal B}_0}\left( \xi \right) \right) =d\left(
\Phi _{{\cal E}_k}\left( \xi \right) \right) d\left( \Phi _{{\cal J}_0}
\left( \xi \right) \right)
$$
have the following form:
$$
d\left( \Phi _{{\cal J}_0}\left( \xi \right) \right) =\exp \{\xi
H_{1+2}\otimes E_{1+2}\},
$$
$$
d\left( \Phi _{{\cal E}_k}\left( \xi \right) \right) =\exp \left\{ \xi
\left( E_1\otimes E_2+2\sum_{j>2}^NE_{1\pm j}\otimes E_{2\mp j}\right)
\right\} .
$$
The $R$ -matrix can be easily obtained using the general formula (\ref
{chrmat}),
$$
\begin{array}{l}
R_{
{\cal B}_0}=I\otimes I-\xi \left( H_{1+2}\wedge E_{1+2}+E_1\otimes
E_2+2\sum_{j>2}^NE_{1\pm j}\otimes E_{2\mp j}\right) + \\ +\frac 12\xi
^2\left( E_1^2\otimes E_2^2+E_2^2\otimes E_1^2+2E_{1+2}\otimes
E_{1+2}-2E_2E_1\otimes E_1E_2\right) + \\
+2\xi ^2\sum_{j>2}^N\left( E_{1\pm j}E_{1\mp j}\otimes E_{2\mp j}E_{2\pm
j}+E_{2\mp j}E_{2\pm j}\otimes E_{1\pm j}E_{1\mp j}\right.  \\
\left. -2E_{2\mp j}E_{1\pm j}\otimes E_{1\pm j}E_{2\mp j}\right)
\end{array}
$$
The corresponding matrix $d(\rho)$ is
$$
d\left( \rho _{{\cal B}_0} \right) = H_{1+2}\otimes E_{1+2}+E_1\otimes
E_2+2\sum_{j>2}^NE_{1\pm j}\otimes E_{2\mp j}.
$$
This expression together with the general formula (\ref{rymatr}) defines
the final result -- the set of $R_{\cal Y}$-matrices for the deformed
Yangians $Y_{{\cal B}_0}(so(2N+1))$.


\section{Conclusions}

We have demonstrated that for some types of simple Lie algebras
the Yangians deformed by chains of extended twists are completely
determined by the matrix $d(\rho)$ which is the logarithm of the twisting
element in the defining representation of an algebra. This gives
rise to a set $R_{\cal Y}$ of solutions to the quantum Yang-Baxter
equation and, therefore, to a series of integrable models with the
corresponding local Hamiltonians.

When the chain has the index $p>1$ the $R_{\cal Y}$-matrices naturally
acquire the parameters $\eta_k$. Moreover, in chains of extended twists
one can cut any number of factors from the left, i.e. use the discrete
parametrization of the last nontrivial sequence ${\cal F}_{{\cal B}_p}$.


\section{Acknowledgements}
I am grateful to professors P.P.Kulish and A.Stolin for
valuable discussions on the subject. It is a pleasure for me to thank
the organizers of the International Seminar "Supersymmetries and Quantum
Symmetries" (SQS'99, 27-31 July, 1999) for warm hospitality.

This work has been partially supported by the Russian Foundation for
Basic Research under the grant N 97-01-01152.


\end{document}